\documentclass{gtart_a}
\pdfoutput=1
\usepackage{pinlabel}

\title{Manifolds with non-stable fundamental groups at infinity, III}

\author{C\,R Guilbault}
\givenname{C\,R}
\surname{Guilbault}
\address{Department of Mathematical Sciences\\
University of Wisconsin-Milwaukee\\\newline
Milwaukee, WI 53201\\USA}
\email{craigg@uwm.edu}
\urladdr{}

\author{F\,C Tinsley}
\givenname{F\,C}
\surname{Tinsley}
\address{Department of Mathematics\\The Colorado College\\\newline
Colorado Springs\\Colorado 80903\\USA}
\email{ftinsley@coloradocollege.edu}
\urladdr{}

\volumenumber{10}
\issuenumber{}
\publicationyear{2006}
\papernumber{14}
\lognumber{0274}
\startpage{541}
\endpage{556}

\doi{}
\MR{}
\Zbl{}

\arxivreference{math.GT/0505517}
\arxivpassword{}   

\keyword{manifold}
\keyword{end}
\keyword{tame}
\keyword{inward tame}
\keyword{open collar}
\keyword{pseudo-collar}
\keyword{semistable}
\keyword{Mittag-Leffler}
\keyword{perfect group}
\keyword{perfectly semistable}
\keyword{Siebenmann's thesis}
\keyword{Wall finiteness obstruction}
\keyword{Quillen's plus construction}
\subject{primary}{msc2000}{57N15}
\subject{primary}{msc2000}{57Q12}
\subject{secondary}{msc2000}{57R65}
\subject{secondary}{msc2000}{57Q10}

\received{24 May 2005}
\revised{}
\accepted{24 March 2006}
\published{27 April 2006}
\publishedonline{27 April 2006}
\proposed{Steven Ferry}
\seconded{Benson Farb, Martin Bridson}
\corresponding{}
\editor{}
\version{}

\arxivreference{math.GT/0505517}
\arxivpassword{}

\AtBeginDocument{\let\bar\wbar}

\makeatletter
\def\cnewtheorem#1[#2]#3{\newtheorem{#1}{#3}[section]
\expandafter\let\csname c@#1\endcsname\c@theorem}

\newtheorem{theorem}{Theorem}[section]
\cnewtheorem{corollary}[theorem]{Corollary}
\cnewtheorem{proposition}[theorem]{Proposition}
\cnewtheorem{lemma}[theorem]{Lemma}

\theoremstyle{remark}

\newtheorem{claim}{Claim}

\cnewtheorem{exercise}[theorem]{Exercise}

\newtheorem{remark}{Remark}

\makeatother

\numberwithin{equation}{section}

\begin{document}

\begin{abstract}
We continue our study of ends non-compact manifolds. The over-arching aim is
to provide an appropriate generalization of Siebenmann's famous collaring
theorem that applies to manifolds having non-stable fundamental group systems
at infinity. In this paper a primary goal is finally achieved; namely, a
complete characterization of pseudo-collarability for manifolds of dimension
at least 6.

\end{abstract}

\maketitle

\section{Introduction}

This is the third in a series of papers aimed at generalizing Siebenmann's
famous PhD thesis \cite{Si} so that the results apply to manifolds with
nonstable fundamental groups at infinity. Siebenmann's work provides necessary
and sufficient conditions for an open manifold of dimension $\geq6$ to contain
an open collar neighborhood of infinity, ie, a manifold neighborhood of
infinity $N$ such that $N\approx\partial N\times\lbrack0,1)$. Clearly, a
stable fundamental group at infinity is necessary in order for such a
neighborhood to exist. Hence, our first task was to identify a useful, but
less rigid, `end structure'\ to aim for. We define a manifold $N^{n}$ with
compact boundary to be a \emph{homotopy collar} provided $\partial
N^{n}\hookrightarrow N^{n}$ is a homotopy equivalence. Then define a
\emph{pseudo-collar} to be a homotopy collar which contains arbitrarily small
homotopy collar neighborhoods of infinity. An open manifold (or more
generally, a manifold with compact boundary) is \emph{pseudo-collarable }if it
contains a pseudo-collar neighborhood of infinity. Obviously, an open collar
is a special case of a pseudo-collar.  Guilbault \cite{Gu2} contains a detailed
discussion of pseudo-collars, including motivation for the definition and a
variety of examples---both pseudo-collarable and non-pseudo-collarable. In
addition, a set of three conditions (see below) necessary for
pseudo-collarability---each analogous to a condition from Siebenmann's
original theorem---was identified there. A primary goal became establishment
of the sufficiency of these conditions. At the time \cite{Gu2} was written, we
were only partly successful at attaining that goal. We obtained an existence
theorem for pseudo-collars, but only by making an additional assumption
regarding the second homotopy group at infinity. In this paper we eliminate
that hypothesis; thereby obtaining the following complete characterization.

\begin{theorem}
[Pseudo-collarability Characterization Theorem]\label{PCT}A one ended
$n$--manifold $M^{n}$ ($n\geq6$) with compact boundary is pseudo-collarable iff
each of the following conditions holds:

\begin{enumerate}
\item $M^{n}$ is inward tame at infinity,

\item $\pi_{1}(\varepsilon(M^{n}))$ is perfectly semistable, and

\item $\sigma_{\infty}\left(  M^{n}\right)  =0\in\widetilde{K}_{0}\left(
\pi_{1}\left(  \varepsilon(M^{n}\right)  \right)  )$.\bigskip
\end{enumerate}
\end{theorem}

\begin{remark}
While it is convenient to focus on one ended manifolds, the above theorem
actually applies to all manifolds with compact boundary. This is true because
an inward tame manifold with compact boundary has only finitely many ends.
(See Section 3 of Guilbault--Tinsley \cite{GuTi}.)  Hence, \fullref{PCT} may be applied to
each end individually. Manifolds with non-compact boundaries are an entirely
different story and will not be discussed here. A detailed discussion of that
situation will be provided in Guilbault \cite{Gu3}.\medskip
\end{remark}

\begin{remark}
\label{r2}A side benefit of our new proof is the inclusion of the $n=6$ case.
In fact, our proof is also valid in dimension five when all of the
groups involved are `good' in the sense of Freedman and Quinn
\cite{FQ}. In that dimension the pseudo-collar structure obtained is
purely topological---as opposed to PL or smooth. This parallels the
status of Siebenmann's theorem in dimension $5$. We discuss dimensions
$\leq4$ at the end of this section.
\end{remark}

The condition of \emph{inward tameness }means that each neighborhood of
infinity can be pulled into a compact subset of itself, or equivalently, that
$M^{n}$ contains arbitrarily small neighborhoods of infinity which are
finitely dominated. Next let $\pi_{1}(\varepsilon\left(  M^{n}\right)  )$
denote the inverse system of fundamental groups of neighborhoods of infinity.
Such a system is \emph{semistable} if it is equivalent to a system in which
all bonding maps are surjections. If, in addition, it can be arranged that the
kernels of these bonding maps are perfect groups, then the system is
\emph{perfectly semistable}. The obstruction $\sigma_{\infty}\left(
M^{n}\right)  \in\widetilde{K}_{0}\left(  \pi_{1}\left(  \varepsilon
(M^{n}\right)  \right)  )$ vanishes precisely when each (clean) neighborhood
of infinity has finite homotopy type. More detailed formulations of these
definitions will be given in \fullref{definitions}.

Conditions (1)--(3) correspond directly to the three conditions identified by
Siebenmann as necessary and sufficient for the existence of an actual open
collar neighborhood of infinity for manifolds of dimension $\geq6$. (His
original version combined the first two into a single assumption.) Indeed,
condition (1) is precisely one of his conditions, condition (2) is a relaxation
of his $\pi_{1}$--stability condition, and condition (3) is the natural
reformulation of his third condition to the situation where $\pi_{1}\left(
\varepsilon(M^{n}\right)  )$ is not necessarily stable.

In the second paper of this series \cite{GuTi}, we focused our attention on
the interdependence of conditions (1)--(3). Specifically, it seemed that
condition (2) might be implied by condition (1), or by a combination of (1) and
(3). This turned out to be partly true. We showed that, for manifolds with
compact boundary, inward tameness implies $\pi_{1}$--semistability. However,
that paper also presents examples satisfying both (1) and (3) which do not have
\emph{perfectly}\thinspace semistable fundamental group at infinity---and thus
are not pseudo-collarable. Those results solidified conditions (1)--(3) as the
best hope for a complete characterization of manifolds with pseudo-collarable ends.

The proof of the Pseudo-collar Characterization Theorem is based upon the
proof of the `Main Existence Theorem' of \cite{Gu2}, which was based on
Siebenmann's original work. The primary task of this paper is to redo the
final step of our earlier proof without assuming $\pi_{2}$--semistability. In
an interesting twist, our new strategy results in a proof that more closely
resembles Siebenmann's original argument than its predecessor. Even so, the
reader would be well served to have a copy of \cite{Gu2} available.

\begin{remark}
When $\pi_{1}(\varepsilon(M^{n}))$ is stable, conditions (1)--(3) become
identical to Siebenmann's conditions. Thus, an application of \cite{Si} tells
us that every pseudo-collar with stable fundamental group at infinity contains
a genuine collar. This fact can also be obtained by a relatively simple direct
argument. Thus, one may view Siebenmann's theorem as a special case of the
Pseudo-collarability Characterization Theorem. For completeness, we have
included that direct argument as \fullref{pseudocollarimpliescollar}
in the following section.
\end{remark}

\begin{remark}
For irreducible $3$--dimensional manifolds with compact boundary the
assumption of inward tameness, by itself, implies the existence of an
open collar neighborhood of infinity, Tucker \cite{Tu}. By contrast,
in dimension $4$, Kwasik and Schultz \cite{KS} have given examples
where Siebenmann's Collaring Theorem fails. Since these examples have
`good'\ fundamental groups at infinity,
\fullref{pseudocollarimpliescollar} (or a quick review of \cite{KS})
shows that \fullref{PCT} also fails in dimension $4$.
\end{remark}

The first author wishes to acknowledge support from NSF Grant DMS-0072786.

\section{Definitions and terminology}\label{definitions}

In this section we briefly review most of the terminology and notation needed
in the remainder of the paper. It is divided into two subsections---the first
devoted to inverse sequences of groups and the second to the topology of ends
of manifolds.

\subsection{Algebra of inverse sequences}

Throughout this section all arrows denote homomorphisms, while arrows of the
type $\twoheadrightarrow$ or $\twoheadleftarrow$ denote surjections. The
symbol $\cong$ denotes isomorphisms.

Let
\[
G_{0}\overset{\lambda_{1}}{\longleftarrow}G_{1}\overset{\lambda_{2}%
}{\longleftarrow}G_{2}\overset{\lambda_{3}}{\longleftarrow}\cdots
\]
be an inverse sequence of groups and homomorphisms. A \emph{subsequence} of
$\left\{  G_{i},\lambda_{i}\right\}  $ is an inverse sequence of the form
\[
G_{i_{0}}\overset{\lambda_{i_{0}+1}\circ\cdots\circ\lambda_{i_{1}}%
}{\longleftarrow}G_{i_{1}}\overset{\lambda_{i_{1}+1}\circ\cdots\circ
\lambda_{i_{2}}}{\longleftarrow}G_{i_{2}}\overset{\lambda_{i_{2}+1}\circ
\cdots\circ\lambda_{i_{3}}}{\longleftarrow}\cdots.
\]
In the future we will denote a composition $\lambda_{i}\circ\cdots\circ
\lambda_{j}$ ($i\leq j$) by $\lambda_{i,j}$.

Sequences $\left\{  G_{i},\lambda_{i}\right\}  $ and $\left\{  H_{i},\mu
_{i}\right\}  $ are \emph{pro-equivalent} if, after passing to subsequences,
there exists a commuting diagram:
\[%
\begin{array}
[c]{ccccccc}%
G_{i_{0}} & \overset{\lambda_{i_{0}+1,i_{1}}}{\longleftarrow} & G_{i_{1}} &
\overset{\lambda_{i_{1}+1,i_{2}}}{\longleftarrow} & G_{i_{2}} & \overset
{\lambda_{i_{2}+1,i_{3}}}{\longleftarrow} & \cdots\\
& \nwarrow\quad\swarrow &  & \nwarrow\quad\swarrow &  & \nwarrow\quad\swarrow
& \\
& H_{j_{0}} & \overset{\mu_{j_{0}+1,j_{1}}}{\longleftarrow} & H_{j_{1}} &
\overset{\mu_{j_{1}+1,j_{2}}}{\longleftarrow} & H_{j_{2}} & \cdots
\end{array}
.
\]
Clearly an inverse sequence is pro-equivalent to any of its subsequences. To
avoid tedious notation, we often do not distinguish $\left\{  G_{i}%
,\lambda_{i}\right\}  $ from its subsequences. Instead we simply assume that
$\left\{  G_{i},\lambda_{i}\right\}  $ has the desired properties of a
preferred subsequence---often prefaced by the words `after passing to a
subsequence and relabelling'.

The \emph{inverse limit }of a sequence $\left\{  G_{i},\lambda_{i}\right\}  $
is a subgroup of $\prod G_{i}$ defined by
\[
\underleftarrow{\lim}\left\{  G_{i},\lambda_{i}\right\}  =\left\{  \left.
\left(  g_{0},g_{1},g_{2},\cdots\right)  \in\prod_{i=0}^{\infty}G_{i}\right|
\lambda_{i}\left(  g_{i}\right)  =g_{i-1}\right\}  .
\]
Notice that for each $i$, there is a \emph{projection homomorphism}
$p_{i}\co \underleftarrow{\lim}\left\{  G_{i},\lambda_{i}\right\}  \rightarrow
G_{i}$. It is a standard fact that pro-equivalent inverse sequences have
isomorphic inverse limits.

An inverse sequence $\left\{  G_{i},\lambda_{i}\right\}  $ is \emph{stable} if
it is pro-equivalent to an inverse sequence $\left\{  H_{i},\mu_{i}\right\}  $
for which each $\mu_{i}$ is an isomorphism. A more usable formulation is that
$\left\{  G_{i},\lambda_{i}\right\}  $ is stable if, after passing to a
subsequence and relabelling, there is a commutative diagram of the form
\begin{equation}%
\begin{array}
[c]{ccccccccc}%
G_{0} & \overset{\lambda_{1}}{\longleftarrow} & G_{1} & \overset{\lambda_{2}%
}{\longleftarrow} & G_{2} & \overset{\lambda_{3}}{\longleftarrow} & G_{3} &
\overset{\lambda_{4}}{\longleftarrow} & \cdots\\
& \nwarrow\quad\swarrow &  & \nwarrow\quad\swarrow &  & \nwarrow\quad\swarrow
&  &  & \\
& im(\lambda_{1}) & \longleftarrow & im(\lambda_{2}) & \longleftarrow &
im(\lambda_{3}) & \longleftarrow & \cdots &
\end{array}
\tag{$\ast$}\label{eq:ast}%
\end{equation}
where each bonding map in the bottom row (obtained by restricting the
corresponding $\lambda_{i}$) is an isomorphism. If $\left\{  H_{i},\mu
_{i}\right\}  $ can be chosen so that each $\mu_{i}$ is an epimorphism, we say
that our inverse sequence is \emph{semistable }(or \emph{Mittag--Leffler,
}or\emph{\ pro-epimorphic}). In this case, it can be arranged that the
restriction maps in the bottom row of \eqref{eq:ast} are epimorphisms. Similarly, if
$\left\{  H_{i},\mu_{i}\right\}  $ can be chosen so that each $\mu_{i}$ is a
monomorphism, we say that our inverse sequence is \emph{pro-monomorphic}; it
can then be arranged that the restriction maps in the bottom row of \eqref{eq:ast}
are monomorphisms. It is easy to see that an inverse sequence that is
semistable and pro-monomorphic is stable.

Recall that a \emph{commutator} element of a group $H$ is an element of the
form $x^{-1}y^{-1}xy$ where $x,y\in H$; and the \emph{commutator subgroup} of
$H,$ denoted $\left[  H,H\right]  $, is the subgroup generated by all of its
commutators. The group $H$ is \emph{perfect} if $\left[  H,H\right]  =H$. An
inverse sequence of groups is \emph{perfectly semistable} if it is
pro-equivalent to an inverse sequence
\[
G_{0}\overset{\lambda_{1}}{\twoheadleftarrow}G_{1}\overset{\lambda_{2}%
}{\twoheadleftarrow}G_{2}\overset{\lambda_{3}}{\twoheadleftarrow}\cdots
\]
of finitely presentable groups and surjections where each $\ker\left(
\lambda_{i}\right)  $ perfect. The following shows that inverse sequences of
this type behave well under passage to subsequences.

\begin{lemma}
\label{L1}A composition of surjective group homomorphisms, each having perfect
kernels, has perfect kernel. Thus, if an inverse sequence of surjective group
homomorphisms has the property that the kernel of each bonding map is perfect,
then each of its subsequences also has that property.
\end{lemma}

\begin{proof}
See  \cite[Lemma 1]{Gu2}.
\end{proof}

\subsection{Topology of ends of manifolds}

Throughout this paper, $\approx$ will represent homeomorphism, while $\simeq$
will indicate homotopic maps or homotopy equivalent spaces. The word
\emph{manifold} means \emph{manifold with (possibly empty) boundary}. A
manifold is \emph{open} if it is non-compact and has no boundary. We will
restrict our attention to manifolds with compact boundaries. This prevents the
`topology at infinity'\ of our manifold from getting entangled with the
topology at infinity of its boundary. Manifolds with noncompact boundaries
will be addressed in \cite{Gu3}.

For convenience, all manifolds are assumed to be PL. Analogous results may be
obtained for smooth or topological manifolds in the usual ways. Occasionally
we will observe that a theorem remains valid in dimension $4$ or $5$. Results
of this sort usually require the purely topological $4$--dimensional techniques
developed by Freedman \cite{FQ}; thus, the corresponding conclusions are
only topological. The main focus of this paper, however, is on dimensions
$\geq6$.

Let $M^{n}$ be a manifold with compact (possibly empty) boundary. A set
$N\subset M^{n}$ is a \emph{neighborhood of infinity} if $\overline{M^{n}-N}$
is compact. A neighborhood of infinity $N$ is \emph{clean} if

\begin{itemize}
\item $N$ is a closed subset of $M^{n}$,

\item $N\cap\partial M^{n}=\emptyset$, and

\item $N$ is a codimension 0 submanifold of $M^{n}$ with bicollared boundary.
\end{itemize}

\noindent It is easy to see that each neighborhood of infinity contains a
clean neighborhood of infinity.

We say that $M^{n}$ \emph{has }$k$ \emph{ends }if it contains a compactum $C$
such that, for every compactum $D$ with $C\subset D$, $M^{n}-D$ has exactly
$k$ unbounded components, ie, $k$ components with noncompact closures. When
$k$ exists, it is uniquely determined; if $k$ does not exist, we say $M^{n}$
has \emph{infinitely many ends}. If $M^{n}$ is $k$--ended, then it contains a
clean neighborhood of infinity $N$ consisting of $k$ connected components,
each of which is a one ended manifold with compact boundary. Thus, when
studying manifolds with finitely many ends, it suffices to understand the
\emph{one ended} situation. That is the case in this paper, where our standard
hypotheses ensure finitely many ends. See \cite[Proposition 3.1]{GuTi}.

A connected clean neighborhood of infinity with connected boundary is called a
\emph{0--neighborhood of infinity}. If $N$ is clean and connected but has more
than one boundary component, we may choose a finite collection of disjoint
properly embedded arcs in $N$ that connect those components. Deleting from $N$
the interiors of regular neighborhoods of these arcs produces a $0$%
--neighborhood of infinity $N_{0}\subset N$.

A nested sequence $N_{0}\supset N_{1}\supset N_{2}\supset\cdots$ of
neighborhoods of infinity is \emph{cofinal }if $\bigcap_{i=0}^{\infty}%
N_{i}=\emptyset$. For any one ended manifold $M^{n}$, one may easily obtain a
cofinal sequence of $0$--neighborhoods of infinity.

We say that $M^{n\text{ }}$is \emph{inward tame }at infinity if, for
arbitrarily small neighborhoods of infinity $N$, there exist homotopies
$H\co N\times\left[  0,1\right]  \rightarrow N$ such that $H_{0}=id_{N}$ and
$\overline{H_{1}\left(  N\right)  }$ is compact. Thus inward tameness means
each neighborhood of infinity can be pulled into a compact subset of itself.

Recall that a space $X$ is \emph{finitely dominated} if there exists a finite
complex $K$ and maps $u\co X\rightarrow K$ and $d\co K\rightarrow X$ such that
$d\circ u\simeq id_{X}$. The following lemma uses this notion to offer
equivalent formulations of `inward tameness'.

\begin{lemma}
{\rm\cite[Lemma 2.4]{GuTi}}\qua For a manifold $M^{n}$, the following are equivalent.

\begin{enumerate}
\item $M^{n}$ is inward tame at infinity.

\item Each clean neighborhood of infinity in $M^{n}$ is finitely dominated.

\item For each cofinal sequence $\left\{  N_{i}\right\}  $ of clean
neighborhoods of infinity, the inverse sequence
\[
N_{0}\overset{j_{1}}{\hookleftarrow}N_{1}\overset{j_{2}}{\hookleftarrow}%
N_{2}\overset{j_{3}}{\hookleftarrow}\cdots
\]
is pro-homotopy equivalent to an inverse sequence of finite polyhedra.
\end{enumerate}
\end{lemma}

Given a nested cofinal sequence $\left\{  N_{i}\right\}  _{i=0}^{\infty}$ of
connected neighborhoods of infinity, base points $p_{i}\in N_{i}$, and paths
$\alpha_{i}\subset N_{i}$ connecting $p_{i}$ to $p_{i+1}$, we obtain an
inverse sequence:
\[
\pi_{1}\left(  N_{0},p_{0}\right)  \overset{\lambda_{1}}{\longleftarrow}%
\pi_{1}\left(  N_{1},p_{1}\right)  \overset{\lambda_{2}}{\longleftarrow}%
\pi_{1}\left(  N_{2},p_{2}\right)  \overset{\lambda_{3}}{\longleftarrow}%
\cdots.
\]
\smallskip Here, each $\lambda_{i+1}\co \pi_{1}\left(  N_{i+1},p_{i+1}\right)
\rightarrow\pi_{1}\left(  N_{i},p_{i}\right)  $ is the homomorphism induced by
inclusion followed by the change of base point isomorphism determined by
$\alpha_{i}$. The obvious singular ray obtained by piecing together the
$\alpha_{i}$'s is often referred to as the \emph{base ray }for the inverse
sequence. Provided the sequence is semistable, one can show that its
pro-equivalence class does not depend on any of the choices made above. We
refer to the pro-equivalence class of this sequence as the \emph{fundamental
group system at infinity} for $M^{n}$ and denote it by $\pi_{1}\left(
\varepsilon\left(  M^{n}\right)  \right)  $. We denote the inverse limit of
this sequence by $\check{\pi}_{1}\left(  \varepsilon\left(  M^{n}\right)
\right)  $. (In the absence of semistability, the pro-equivalence class of the
inverse sequence depends on the choice of base ray, and hence, this choice
becomes part of the data.) It is easy to see how the same procedure may also
be used to define $\pi_{k}\left(  \varepsilon\left(  M^{n}\right)  \right)  $
and $\check{\pi}_{k}\left(  \varepsilon\left(  M^{n}\right)  \right)  $for
$k>1$.

In \cite{Wa}, Wall shows that each finitely dominated connected space $X$
determines a well-defined element $\sigma\left(  X\right)  $ lying in
$\widetilde{K}_{0}\left(  \mathbb{Z}\left[  \pi_{1}X\right]  \right)  $ (the
group of stable equivalence classes of finitely generated projective
$\mathbb{Z}\left[  \pi_{1}X\right]  $--modules under the operation induced by
direct sum) that vanishes if and only if $X$ has the homotopy type of a finite
complex. Given a nested cofinal sequence $\left\{  N_{i}\right\}
_{i=0}^{\infty}$ of connected clean neighborhoods of infinity in an inward
tame manifold $M^{n}$, we have a Wall obstruction $\sigma(N_{i})$ for each
$i$. These may be combined into a single obstruction
\begin{align*}
\sigma_{\infty}(M^{n})&=\left(  -1\right)  ^{n}(\sigma(N_{0}),\sigma
(N_{1}),\sigma(N_{2}),\cdots)\\&\in\widetilde{K}_{0}\left(  \pi_{1}\left(
\varepsilon\left(  M^{n}\right)  \right)  \right)  \equiv\underleftarrow{\lim
}\widetilde{K}_{0}\left(  \mathbb{Z}\left[  \pi_{1}N_{i}\right]  \right)
\end{align*}
that is well-defined and which vanishes if and only if each clean
neighborhood of infinity in $M^{n}$ has finite homotopy type. See
Chapman and Siebenmann \cite{CS} for details.

We conclude this section by providing a direct proof of the following:

\begin{proposition}
\label{pseudocollarimpliescollar}Every pseudo-collar of dimension $\geq5$ with
stable fundamental group contains an open collar neighborhood of infinity. In
dimension $4$ this remains true in the topological category provided the
fundamental group at infinity is good.
\end{proposition}

\begin{proof}
If $N^{n}$ is a connected pseudo-collar, then by definition there exists a
cofinal sequence $N^{n}=N_{1}\supseteq N_{2}\supseteq N_{3}\supseteq\cdots$ of
homotopy collar neighborhoods of infinity. Letting $W_{i}=N_{i}-\overset
{\circ}{N}_{i+1}$ for each $i$, divides $N^{n}$ into a countable sequence
$\left\{  \left(  W_{i},\partial N_{i},\partial N_{i+1}\right)  \right\}  $ of
cobordisms with the property that each $\partial N_{i}\hookrightarrow W_{i}$
is a homotopy equivalence. (We call these \emph{one-sided} $h$--cobordisms.)
Although $\partial N_{i+1}\hookrightarrow W_{i}$ needn't be a homotopy
equivalence, an argument involving duality in the universal cover of $W_{i}$
(see \cite[Theorem 2.5]{GuTi}), implies that $\pi_{1}\left(  \partial
N_{i+1}\right)  \rightarrow$ $\pi_{1}\left(  W_{i}\right)  $ is surjective for
each $i$. By commutativity of%
\[%
\begin{array}
[c]{ccc}%
\pi_{1}\left(  \partial N_{i+1}\right)   & \twoheadrightarrow & \pi_{1}\left(
W_{i}\right)  \\
\cong\downarrow &  & \cong\downarrow\\
\pi_{1}\left(  N_{i+1}\right)   & \longrightarrow & \pi_{1}\left(
N_{i}\right)
\end{array}
\]
each bonding map in the following sequence is surjective
\[
\pi_{1}\left(  N_{1}\right)  \twoheadleftarrow\pi_{1}\left(  N_{2}\right)
\twoheadleftarrow\pi_{1}\left(  N_{3}\right)  \twoheadleftarrow\cdots.
\]
The only way that an inverse sequence of surjections can be stable is
that eventually all bonding homomorphisms are isomorphisms. Choose
$i_{0}$ sufficiently large that $\pi_{1}\left( N_{i+1}\right)
\rightarrow\pi _{1}\left( N_{i}\right) $ is an isomorphism for all
$i\geq i_{0}$. Then for each $i\geq i_{0}$, $\left( W_{i},\partial
N_{i},\partial N_{i+1}\right) $ is a genuine $h$--cobordism. If each is
a product, we may piece the product structures together to obtain an
open collar structure on $N_{i_{0}}$.  Otherwise we will apply the
`weak $h$--cobordism theorem' (Stallings \cite{St} or Connell
\cite{Con}) to rechoose the cobordisms so that each is a product. To
accomplish this, assume that $n\geq5$. Then, by the weak $h$--cobordism
theorem, there is a homeomorphism $h\co \partial N_{i_{0}}\times\lbrack
0,1)\rightarrow W_{i_{0}}-\partial N_{i_{0}+1}$. Choose $t$ close to
$1$; then replace $N_{i_{0}+1}$ with
$N_{i_{0}+1}^{\prime}=N_{i_{0}}-h\left( \partial
N_{i_{0}}\times\lbrack0,t)\right) $ and $W_{i_{0}}$ with
$W_{i_{0}}^{\prime
}=N_{i_{0}}-\overset{\circ}{N^{\prime}}_{i_{0}+1}=h\left( \partial
N_{i}\times\lbrack0,t]\right) $. Next apply the same procedure to
$N_{i_{0}+1}^{\prime}$, $N_{i_{0}+2}$ and the $h$--cobordism between
$\partial N_{i_{0}+1}^{\prime}$ and $\partial N_{i_{0}+2}$ to rechoose
that cobordism so that it is a product. Continue this procedure
inductively to arrange that all of the $h$--cobordisms are products.
If $n=4$, the weak $h$--cobordism theorem (and hence, the above proof)
is still valid provided the common fundamental group is good and the
conclusion we seek is only topological \cite[Corollary 3.5]{Gu1}.
\end{proof}

\section{One sided $h$--cobordisms and the Plus Construction}

As noted above, a pseudo-collar structure on the end of a manifold allows one
to express that end as a countable union of compact `one-sided $h$--cobordisms'
$\left\{  \left(  W_{i},A_{i},B_{i}\right)  \right\}  _{i=1}^{\infty}$ in the
sense that each inclusion $A_{i}\hookrightarrow W_{i}$ is a homotopy
equivalence. For any such cobordism, $\pi_{1}\left(  B_{i}\right)
\rightarrow\pi_{1}\left(  W_{i}\right)  $ is surjective and has perfect kernel
(again see \cite[Theorem 2.5]{GuTi}).  Quillen's famous `plus construction'
(\cite{Qu} or \cite[Section 11.1]{FQ}) provides a partial converse to that observation.

\begin{theorem}
[The Plus Construction]\label{PC}Let $B$ be a closed $\left(  n-1\right)
$--manifold $\left(  n\geq6\right)  $ and $h\co \pi_{1}\left(  B\right)
\rightarrow H$ a surjective homomorphism onto a finitely presented group such
that $\ker\left(  h\right)  $ is perfect. Then there exists a compact
$n$--dimensional cobordism $\left(  W,A,B\right)  $ such that $\ker\left(
\pi_{1}\left(  B\right)  \rightarrow\pi_{1}\left(  W\right)  \right)  =\ker
h$, and $A\hookrightarrow W$ is a simple homotopy equivalence. These
properties determine $W$ uniquely up to homeomorphism rel $B$. If $n=5$, this
result is still valid (in the topological category) provided the group $H$ is good.
\end{theorem}

When a one-sided $h$--cobordism has trivial Whitehead torsion, ie, when the
corresponding homotopy equivalence is simple, we refer to it as a \emph{plus
cobordism}.

In \cite{Gu2} techniques borrowed from the proof of \fullref{PC} were used
in obtaining pseudo-collar structures. In the current paper we apply the plus
construction itself to isolate the difficulties caused by non-stability of the
fundamental group. The key technical tool is the following:

\begin{theorem}
[The Embedded Plus Construction]\label{EPC}Let $R$ be a connected manifold of
dimension $\geq6$; $B$ be a closed component of $\partial R$; and
\[
G\subseteq\ker\left(  \pi_{1}\left(  B\right)  \rightarrow\pi_{1}\left(
R\right)  \right)
\]
a perfect group which is the normal closure in $\pi_{1}\left(  B\right)  $ of
a finite set of elements. Then there exists a plus cobordism $\left(
W,A,B\right)  $ embedded in $R$ which is the identity on $B$ for which
$\ker\left(  \pi_{1}\left(  B\right)  \rightarrow\pi_{1}\left(  W\right)
\right)  =G$. If $n=5$ and $\pi_{1}\left(  B\right)  /G$ is good, the
conclusion is still valid provided we work in the topological category.
\end{theorem}

\begin{proof}
Let $\left(  W,A,B\right)  $ be the plus cobordism promised by \fullref{PC} for the homomorphism $\pi_{1}\left(  B\right)  \rightarrow\pi
_{1}\left(  W\right)  $. We will show (indirectly) how to embed this cobordism
into $R$ in the appropriate manner.
Let $R_{1}=R\cup_{B}W$, the manifold obtained by gluing a copy of $W$ to $R$
along a common copy of $B$.
\begin{claim}
\label{c1}$R\hookrightarrow R_{1}$ is a simple homotopy equivalence.
\end{claim}
Since $A\hookrightarrow W$ is a simple homotopy equivalence, duality implies
that $B\hookrightarrow W$ is a $\mathbb{Z}[\pi_{1}W]$--homology equivalence; in
other words, $H_{\ast}\left(  \widetilde{W},\widehat{B}\right)  =0$, where
$\widehat{B}$ is the preimage of $B$ under the universal covering projection
$\widetilde{W}\rightarrow W$. By Van Kampen's theorem, $R\hookrightarrow
R_{1}$ induces a $\pi_{1}$--isomorphism, so $\widetilde{R}$ may be viewed as a
subset of $\widetilde{R}_{1}.$ So excision implies that $H_{\ast}\left(
\widetilde{R}_{1},\widetilde{R}\right)  \cong H_{\ast}\left(  \widetilde
{W},\widehat{B}\right)  =0$. From there we see that $\pi_{j}\left(
R_{1},R\right)  $ is trivial for all $j$, so $R\hookrightarrow R_{1}$ is a
homotopy equivalence.

The Duality Theorem for Whitehead torsion (Milnor \cite[page 394]{Mi})
may be used to check that this homotopy equivalence is simple. In
particular, even though $B\hookrightarrow W$ is not a homotopy
equivalence, it is a $\mathbb{Z}[\pi_{1}W]$--homology equivalence,
and thus determines an element of $Wh\left( \pi_{1}W\right) $. But
$A\hookrightarrow W$ is a simple homotopy equivalence, so by duality,
both of these inclusions have trivial torsion. By an application of
the Sum Theorem for Whitehead torsion (Cohen \cite[Section 23]{Coh}),
$R\hookrightarrow R_{1}$ also determines the trivial element of
$Wh\left( \pi_{1}R\right) $, so we have the desired simple homotopy
equivalence.  Let
\begin{align*}
X^{n+1}  &  =\left(  R\times\left[  0,\frac{1}{2}\right]  \right)
\cup_{R\times\left\{  \frac{1}{2}\right\}  }\left(  R_{0}\times\left[
\frac{1}{2},1\right]  \right)  ,\smallskip\\
R_{0}  &  =R\times\left\{  0\right\}  \text{, and}\\
Q  &  =(R_{0}\times\left\{  1\right\}  )\cup\left(  A\times\left[  \frac{1}%
{2},1\right]  \right)  \cup\left(  W\times\left\{  \frac{1}{2}\right\}
\right)
\end{align*}
See \fullref{fig1}.%
\begin{figure}[ht!]
\begin{center}
\labellist
\pinlabel $R$ [b] at 288 665
\pinlabel $W$ [b] at 396 665
\pinlabel $1$ [r] at 182 646
\pinlabel $\frac12$ [r] at 182 575
\pinlabel $0$ [r] at 182 502
\pinlabel $Q$ [l] at 479 632
\pinlabel $R_0$ [b] at 402 523
\pinlabel $X^{n+1}$ at 288 575
\endlabellist
\includegraphics[width=3in]{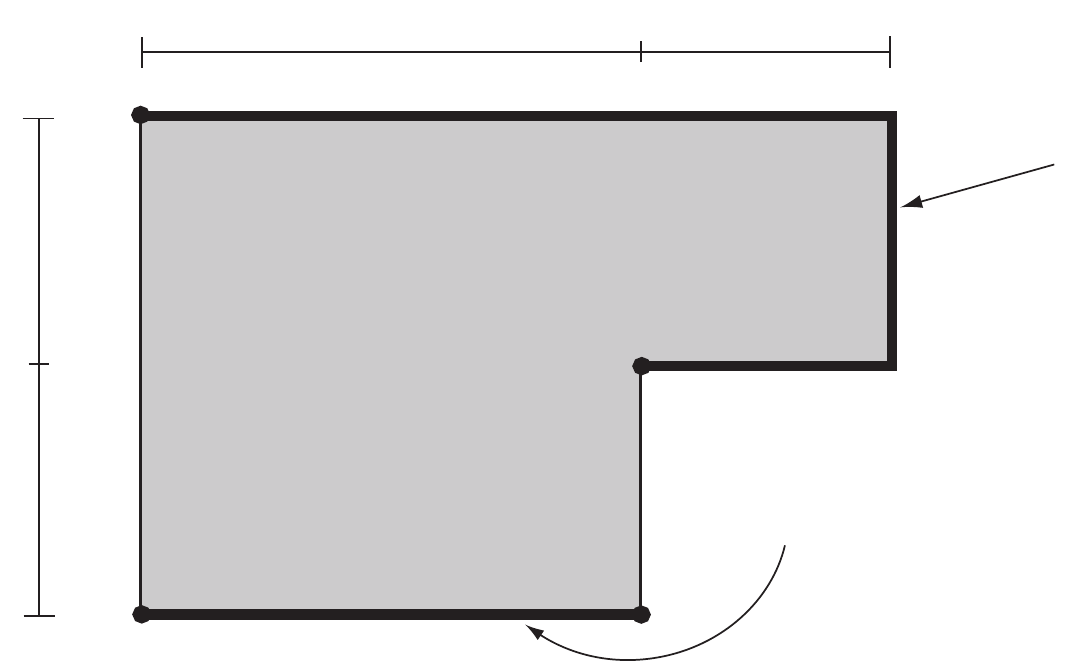}
\caption{The cobordism $\left(  X^{n+1},R_{0},Q\right) $}
\label{fig1}
\end{center}
\end{figure}
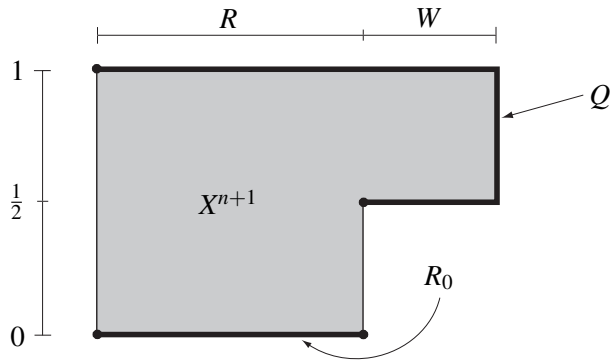
\begin{claim}
\label{c2}$\left(  X^{n+1},R_{0},Q\right)  $ is an $s$--cobordism.
\end{claim}
First note that $\left(  X^{n+1},R_{0}\right)  \simeq\left(  R_{1},R\right)
$, so by Claim \ref{c1}, $R_{0}\hookrightarrow X^{n+1}$ is a homotopy
equivalence. Next observe that $R_{0}\times\left\{  1\right\}  \hookrightarrow
X^{n+1}$ and $R_{0}\times\left\{  1\right\}  \hookrightarrow Q$ are both
homotopy equivalences. The first of these is obvious, while the second follows
from the fact that $A\times\left\{  \frac{1}{2}\right\}  \hookrightarrow
W\times\left\{  \frac{1}{2}\right\}  $ is a homotopy equivalence. It follows
that $Q\hookrightarrow X^{n+1}$ is a homotopy equivalence.

To show that $R_{0}\hookrightarrow X^{n+1}$ has trivial torsion, factor the
inclusion map as $R_{0}\hookrightarrow R\times\lbrack0,1]\hookrightarrow
X^{n+1}$. The first of these inclusions is obviously a simple homotopy
equivalence, and the second is a simple homotopy equivalence by an easy
application of Claim \ref{c1}; thus Claim \ref{c2} follows.

By the relative $s$--cobordism theorem (Rourke and Sanderson
\cite[Chapter 6]{RS}), $X^{n+1}$ is a product, so there is a
homeomorphism (rel boundary) from $Q$ onto $R_{0}$. The image of
$W\times\left\{ \frac{1}{2}\right\} $ under this homeomorphism
provides the desired plus cobordism $\left( W,A,B\right) $ embedded in
$R_{0}$.
\end{proof}

\section[Proof of \ref{PCT}]{Proof of \fullref{PCT}}

We now move to the proof of our main theorem. For a full understanding, the
reader should be familiar with the proof of the Main Existence Theorem of
\cite{Gu2} up to the last few pages---which our current argument will replace.
Those familiar with \cite{Si} will understand the key points. We begin with a
brief review.

Start by assuming only that $M^{n}$ is a one ended manifold with compact
boundary, and that $n$ is at least $5$. (In \cite{Gu2} we took the traditional
route and assumed $M^{n}$ was an open manifold; but this is unnecessary as
long as $\partial M^{n}$ is compact.) Recall that a $0$--neighborhood of
infinity is a \emph{generalized }$1$\emph{--neighborhood of infinity} provided
$\pi_{1}\left(  \partial U\right)  \rightarrow\pi_{1}\left(  U\right)  $ is an
isomorphism. If, in addition, $\pi_{i}\left(  U,\partial U\right)  =0$ for all
$i\leq k$, then $U$ is a \emph{generalized }$k$\emph{--neighborhood of
infinity}.

By the Generalized $\left( n-3\right) $--neighborhoods Theorem
(\cite[Theorem 5] {Gu2}), inward tameness alone allows us to obtain a
cofinal sequence $\left\{ U_{i}\right\} $ of generalized $\left(
n-3\right) $--neighborhoods of infinity in $M^{n}$. Since
\cite[Theorem 1.2]{GuTi} assures that $\pi_{1}\left( \varepsilon\left(
M^{n}\right) \right) $ is semistable, we may also arrange that
$\pi_{1}\left( U_{i}\right) \leftarrow\pi_{1}\left( U_{i+1}\right) $
is surjective for all $i\geq1$. For each $i$ let
$R_{i}=U_{i}-\overset{\circ }{U}_{i+1}$ and consider the collection of
cobordisms $\{\left( R_{i},\partial U_{i},\partial U_{i+1}\right)
\}$. Then

\begin{enumerate}
\item[(i)] Each inclusion $\partial U_{i}\hookrightarrow R_{i}\hookrightarrow
U_{i}$ induces a $\pi_{1}$--isomorphism,

\item[(ii)] $\partial U_{i+1}\hookrightarrow R_{i}$ induces a $\pi_{1}%
$--epimorphism for each $i$,

\item[(iii)] $\pi_{k}(R_{i},\partial U_{i})=0$ for all $k<n-3$ and all $i$, and

\item[(iv)] Each $\left(  R_{i},\partial U_{i},\partial U_{i+1}\right)  $
admits a handle decomposition based on $\partial U_{i}$ containing handles
only of index $\left(  n-3\right)  $ and $\left(  n-2\right)  $.
\end{enumerate}

\noindent The first two observations follow easily from Van Kampen's Theorem.
The third is obtained inductively. First note that by the Hurewicz Theorem
$\pi_{k}(R_{i},\partial U_{i})\cong\pi_{k}(\widetilde{R}_{i},\partial
\widetilde{U}_{i})\cong H_{k}(\widetilde{R}_{i},\partial\widetilde{U}_{i})$,
provided that $\pi_{j}(R_{i},\partial U_{i})$ is trivial for $j<k$. Then
examine the homology long exact sequence for the triple $\left(  \widetilde
{U}_{i},\widetilde{R}_{i},\partial\widetilde{U}_{i}\right)  $ to obtain the
desired result. See \cite[page 561]{Gu2} for details. The fourth observation is
obtained from standard handle theoretic techniques (see \cite{RS}). In
particular, (iii) allows us to eliminate all handles of index $\leq n-4$; then
observation (ii) allows us to eliminate $0$-- and $1$--handles from the
corresponding dual handle decomposition of $R_{i}$ based on $\partial U_{i+1}$.

By observation (iv), each $U_{i}$ admits an infinite handle decomposition
having handles only of index $\left(  n-3\right)  $ and $\left(  n-2\right)
$. Thus, $\left(  U_{i},\partial U_{i}\right)  $ has the homotopy type of a
relative CW pair $\left(  K_{i},\partial U_{i}\right)  $ with $\dim\left(
K_{i}-\partial U_{i}\right)  \leq n-2$. Therefore, if one of the $U_{i}$ is a
generalized $\left(  n-2\right)  $--neighborhood of infinity, then it is a
homotopy collar. Thus, our goal is to improve arbitrarily small $U_{i}$ to
generalized $\left(  n-2\right)  $--neighborhoods of infinity. This must be
done in the above context---in particular, condition (ii) must be preserved. We
accomplish this by altering the $U_{i}$ without changing their fundamental groups.

The next key observation is that, for each $i$, $\pi_{n-2}\left(
U_{i},\partial U_{i}\right)  \cong H_{n-2}(\widetilde{U}_{i},\partial
\widetilde{U}_{i})$\break is a finitely generated projective $\mathbb{Z}[\pi%
_{1}U_{i}]$--module. Moreover, as an element of\break $\widetilde{K}_{0}\left(
\mathbb{Z}[\pi_{1}U_{i}]\right)  $, $\left[  H_{n-2}(\widetilde{U}%
_{i},\partial\widetilde{U}_{i})\right]  =\left(  -1\right)  ^{n}\sigma\left(
U_{i}\right)  $, where $\sigma\left(  U_{i}\right)  $ is the Wall finiteness
obstruction for $U_{i}$. This is the content of \cite[Lemma 13]{Gu2}. As
discussed in \fullref{definitions}, these elements of $\widetilde{K}%
_{0}\left(  \mathbb{Z}[\pi_{1}U_{i}]\right)  $ determine the obstruction
$\sigma_{\infty}\left(  M^{n}\right)  $ found in \fullref{PCT}. By
assuming that $\sigma_{\infty}\left(  M^{n}\right)  $ vanishes, we are given
that each $H_{n-2}(\widetilde{U}_{i},\partial\widetilde{U}_{i})$ is a stably
free $\mathbb{Z}[\pi_{1}U_{i}]$--module. By carving out finitely many trivial
$\left(  n-3\right)  $--handles from each $U_{i}$ we can arrange that these
homology groups are finitely generated free $\mathbb{Z}[\pi_{1}U_{i}%
]$--modules. This can be done so that each remains a generalized $\left(
n-3\right)  $--neighborhood of infinity, and so that none of the fundamental
groups of the neighborhoods of infinity or their boundaries are changed. To
save on notation, we continue to denote this improved collection by $\left\{
U_{i}\right\}  $. See \cite[Lemma 14]{Gu2} for details.

By the finite generation of $H_{n-2}(\widetilde{U}_{i},\partial\widetilde
{U}_{i})$, we may assume (after passing to a subsequence of $\left\{
U_{i}\right\}  $ and relabeling) that $H_{n-2}(\widetilde{R}_{i}%
,\partial\widetilde{U}_{i})\rightarrow H_{n-2}(\widetilde{U}_{i}%
,\partial\widetilde{U}_{i})$ is surjective for each $i$. From there the long
exact sequence for the triple $\left(  \widetilde{U}_{i},\widetilde{R}%
_{i},\partial\widetilde{U}_{i}\right)  $ shows that these surjections are, in
fact, isomorphisms. As above, we may choose a handle decompositions for the
$R_{i}$ based on $\partial U_{i}$ having handles only of index $n-3$ and $n-2$.

From now on, let $i$ be fixed. After introducing some trivial $\left(
n-3,n-2\right)  $--handle pairs, an algebraic lemma and some handle slides
allows us to obtain a handle decomposition of $R_{i}$ based on $\partial
U_{i}$ with $\left(  n-2\right)  $--handles $h_{1}^{n-2},h_{2}^{n-2}%
,\allowbreak\cdots,h_{r}^{n-2}$ and an integer $s\leq r$, such that the
subcollection $\{h_{1}^{n-2},h_{2}^{n-2},\allowbreak\cdots,h_{s}^{n-2}\}$ is a
free $\mathbb{Z}\left[  \pi_{1}R_{i}\right]  $--basis for $H_{n-2}\left(
\widetilde{R}_{i},\partial\widetilde{U}_{i}\right)  $. Then the
corresponding\emph{ }$\mathbb{Z}\left[  \pi_{1}R_{i}\right]  $--cellular chain
complex for $\left(  R_{i},\partial U_{i}\right)  $ may be expressed as%
\begin{equation}
0\rightarrow\left\langle h_{1}^{n-2},\cdots,h_{s}^{n-2}\right\rangle
\oplus\left\langle h_{s+1}^{n-2},\cdots,h_{r}^{n-2}\right\rangle
\overset{\partial}{\longrightarrow}\left\langle h_{1}^{n-3},\cdots,h_{t}%
^{n-3}\right\rangle \rightarrow0 \tag{\dag}\label{eq:dag}%
\end{equation}
where $\left\langle h_{1}^{n-2},\cdots,h_{s}^{n-2}\right\rangle $ represents
the free $\mathbb{Z}\left[  \pi_{1}R_{i}\right]  $--submodule of $\widetilde
{C}_{n-2}$ generated by the corresponding handles; $\left\langle h_{s+1}%
^{n-2},\cdots,h_{r}^{n-2}\right\rangle $ represents the free submodule of
$\widetilde{C}_{n-2}$ generated by the remaining $\left(  n-2\right)
$--handles in $R_{i}$; and
\[
\left\langle h_{1}^{n-3},\cdots,h_{t}^{n-3}\right\rangle =\widetilde{C}_{n-3}%
\]
is the free module generated by the $\left(  n-3\right)  $--handles in $R_{i}$.
Moreover,
\[
H_{n-2}(\widetilde{R}_{i},\partial\widetilde{U}_{i})=\ker\left(
\partial\right)  =\left\langle h_{1}^{n-2},\cdots,h_{s}^{n-2}\right\rangle
\oplus\left\{  0\right\}  ,
\]
and $\partial$ takes $\left\{  0\right\}  \oplus\left\langle h_{s+1}%
^{n-2},\cdots,h_{r}^{n-2}\right\rangle $ injectively into $\left\langle
h_{1}^{n-3},\cdots,h_{t}^{n-3}\right\rangle $. This is the content of Lemma 15
and the following paragraph in \cite{Gu2}.

At this point, we would like to use the fact that $\partial h_{j}^{n-2}=0$ for
each $j=1,\cdots,s$ to slide these handles off all of the $\left(  n-3\right)
$--handles. This would be done by repeated use of the Whitney Lemma in
$\partial_{+}(S\cup h_{1}^{n-3}\cup\cdots\cup h_{t}^{n-3})$ to remove the
collection of attaching spheres $\{\alpha_{j}^{n-3}\}_{j=1}^{s}$ from the belt
spheres $\left\{  \beta_{j}^{2}\right\}  _{j=1}^{t}$ of the $\left(
n-3\right)  $--handles. (Here, $S$ is a closed collar neighborhood of $\partial
U_{i}$ and $\partial_{+}$ indicates the right-hand boundary.) After that, we
would `carve out' these $\left(  n-2\right)  $--handles---those generating the
unwanted $\left(  n-2\right)  $--dimensional homology---in an attempt to obtain
a generalized $\left(  n-2\right)  $--neighborhood of infinity. (This process
will be discussed in detail later.) Unfortunately, the desired application of
the Whitney Lemma is only assured if the collection $\{\alpha_{j}%
^{n-3}\}_{j=1}^{s}$ is $\pi_{1}$--negligible in $\partial_{+}(S\cup h_{1}%
^{n-3}\cup\cdots\cup h_{t}^{n-3})$, ie, the inclusion
\[
\partial_{+}(S\cup h_{1}^{n-3}\cup\cdots\cup h_{t}^{n-3})-\cup_{j=1}^{s}%
\alpha_{j}^{n-3}\hookrightarrow\partial_{+}(S\cup h_{1}^{n-3}\cup\cdots\cup
h_{t}^{n-3})
\]
induces a $\pi_{1}$--isomorphism (see \cite[page 72]{RS}). Moreover, even if the
handles that are generating the unwanted homology can be made to miss the
$\left(  n-3\right)  $--handles, we still must be sure that carving out these
$\left(  n-2\right)  $--handles does not change the fundamental group of our
neighborhood of infinity. Otherwise we will have arranged that the relative
$\mathbb{Z}\left[  \pi_{1}U_{i}\right]  $--homology of our new neighborhood of
infinity is trivial, but $\pi_{1}\left(  U_{i}\right)  $ will be the wrong group.

The above two difficulties are related. To avoid them entirely, we would need
to know that, in the corresponding dual handle decomposition, the $2$--handles
dual to $h_{1}^{n-2},\cdots,h_{s}^{n-2}$ do not kill any non-trivial loops
when they are attached to $\partial U_{i+1}$. Since $\pi_{1}\left(  \partial
U_{i+1}\right)  \rightarrow\pi_{1}\left(  R_{i}\right)  $ is not injective,
that scenario seems highly unlikely.

\begin{remark}
Examples constructed in \cite{GuTi} show that the above problems can indeed occur.
\end{remark}

At this point we begin utilizing Condition (2) of \fullref{PCT}. According
to \cite[Theorem 5]{Gu2}, the collection $\left\{  U_{i}\right\}  $ may then be
chosen so that each homomorphism $\pi_{1}(U_{i+1})\overset{\lambda_{i+1}%
}{\longrightarrow}\pi_{1}(U_{i})$ is surjective with perfect kernel. As noted
earlier, the inclusions $\partial U_{i}\hookrightarrow R_{i}\hookrightarrow
U_{i}$ each induce $\pi_{1}$--isomorphisms. By similar reasoning $\pi
_{1}(\partial U_{i+1})\rightarrow\pi_{1}(R_{i})$ is surjective with the same
kernel as $\lambda_{i+1}$. Call this kernel $K_{i+1}$. By a basic theorem from
combinatorial group theory (see \cite{Si} or \cite[Lemma 3]{Gu2}) $K_{i+1}$ is
the normal closure of a finite collection of elements of $\pi_{1}(\partial
U_{i+1})$. Thus we may apply \fullref{EPC} to $\left(  R_{i},\partial
U_{i},\partial U_{i+1}\right)  $ to obtain a plus cobordism $\left(
W_{i},A_{i},\partial U_{i+1}\right)  $ embedded in $R_{i}$ which is the
identity on $\partial U_{i+1}$ and for which $\ker\left(  \pi_{1}\left(
\partial U_{i+1}\right)  \rightarrow\pi_{1}\left(  W_{i}\right)  \right)
=K_{i+1}$. It follows that $\pi_{1}\left(  W_{i}\right)  \overset{\cong
}{\longrightarrow}\pi_{1}\left(  R_{i}\right)  $

Let $R_{i}^{\prime}=\overline{R_{i}-W_{i}}$. Since $W_{i}$ strong deformation
retracts onto $A_{i}$ we have

\begin{itemize}
\item $\left(  R_{i}^{\prime},\partial U_{i}\right)  \hookrightarrow\left(
R_{i},\partial U_{i}\right)  $ is a homotopy equivalence of pairs, and

\item $\pi_{1}\left(  A_{i}\right)  \overset{\cong}{\longrightarrow}\pi
_{1}\left(  R_{i}^{\prime}\right)  $.
\end{itemize}

\noindent The first property ensures that the inclusion induced maps
\[
H_{n-2}(\widetilde{R}_{i}^{\prime},\partial\widetilde{U}_{i})\rightarrow
H_{n-2}(\widetilde{R}_{i},\partial\widetilde{U}_{i})\rightarrow H_{n-2}%
(\widetilde{U}_{i},\partial\widetilde{U}_{i})
\]
are all isomorphisms; thus, $H_{n-2}(\widetilde{R}_{i}^{\prime},\partial
\widetilde{U}_{i})$ is a free $\mathbb{Z}[\pi_{1}U_{i}]$--module which carries
the $\mathbb{Z}[\pi_{1}U_{i}]$--homology of $\left(  U_{i},\partial
U_{i}\right)  $. By performing the same procedures on the cobordism $\left(
R_{i}^{\prime},\partial U_{i},A\right)  $ as we did earlier on $\left(
R_{i},\partial U_{i},\partial U_{i+1}\right)  $ we may obtain a handle
decomposition of $R_{i}^{\prime}$ based on $\partial U_{i}$ which has handles
only of index $n-3$ and $n-2$. Moreover, we may arrange that the corresponding
cellular chain complex is of the form \eqref{eq:dag} (although the precise numbers of
handles may have changed). We adopt that notation without changing the names
of the handles.

The second property ensures that, under the dual handle decomposition of
$R_{i}^{\prime}$, the $2$--handles dual to $h_{1}^{n-2},\cdots,h_{s}^{n-2}$ do
not kill any non-trivial loops when they are attached to $A$. This means that
the attaching $\left(  n-3\right)  $--spheres of $h_{1}^{n-2},\cdots
,h_{s}^{n-2}$ are all $\pi_{1}$--negligible in $\partial_{+}\left(  S\cup
h_{1}^{n-3}\cup\cdots\cup h_{t}^{n-3}\right)  $. The non-simply connected
Whitney Lemma \cite[page 72]{RS} may now be applied to isotope the attaching
spheres of $h_{1}^{n-2},\cdots,h_{s}^{n-2}$ off all of the belt spheres of the
$\left(  n-3\right)  $--handles. Thus, we may assume that $h_{1}^{n-2}%
,\cdots,h_{s}^{n-2}$ are attached directly to $S$. Let $Q=S\cup\left(
h_{1}^{n-2}\cup\cdots\cup h_{s}^{n-2}\right)  $ and let $V_{i}=\overline
{U_{i}-Q}$. We will show that $V_{i}$ is the desired generalized $\left(
n-2\right)  $--neighborhood of infinity. The first issue involves the
fundamental group. We wish to observe that the fundamental group has not
changed, ie, that $V_{i}\hookrightarrow U_{i}$ induces a $\pi_{1}%
$--isomorphism and that $V_{i}$ is a generalized $1$--neighborhood of infinity.

First, note that $\overline{R_{i}^{\prime}-Q}$ may be obtained from $\partial
V_{i}$ by attaching $(n-3)$-- and $(n-2)$--handles; in particular $\{h_{1}%
^{n-3},\cdots,h_{t}^{n-3}\}$ and $\{h_{s+1}^{n-2},\cdots,h_{r}^{n-2}\}$. Since
$n\geq6$, $\partial V_{i}\hookrightarrow\overline{R_{i}^{\prime}-Q}$ induces a
$\pi_{1}$--isomorphism. By inverting this handle decomposition, it is clear
that $A\hookrightarrow\overline{R_{i}^{\prime}-Q}$ induces a $\pi_{1}%
$--surjection. Moreover, since the composition $A\hookrightarrow\overline
{R_{i}^{\prime}-Q}\hookrightarrow R_{i}^{\prime}$ induces a $\pi_{1}%
$--isomorphism, we also have injectivity; so $A\hookrightarrow\overline
{R_{i}^{\prime}-Q}$ induces a $\pi_{1}$--isomorphism. This also implies that
$\overline{R_{i}^{\prime}-Q}\hookrightarrow R_{i}^{\prime}$ induces a $\pi
_{1}$--isomorphism.

Since $U_{i+1}$ is a generalized $1$--neighborhood of infinity, the Van Kampen
theorem assures us that $A\hookrightarrow W\cup U_{i+1}$ induces a $\pi_{1}%
$--isomorphism. Similar arguments then provide the necessary isomorphisms
$\pi_{1}\left(  \partial V_{i}\right)  \overset{\cong}{\longrightarrow}\pi
_{1}\left(  V_{i}\right)  $ and $\pi_{1}\left(  \partial V_{i}\right)
\overset{\cong}{\longrightarrow}\pi_{1}\left(  U_{i}\right)  $.

Lastly, we verify that $V_{i}$ is a generalized $\left(  n-2\right)
$--neighborhood of infinity. Begin with the long exact sequence for the triple
$\left(  \widetilde{U}_{i},\widetilde{Q},\partial\widetilde{U}_{i}\right)  $.
\[
\cdots\rightarrow H_{k}\left(  \widetilde{Q},\partial\widetilde{U}_{i}\right)
\rightarrow H_{k}\left(  \widetilde{U}_{i},\partial\widetilde{U}_{i}\right)
\rightarrow H_{k}\left(  \widetilde{U}_{i},\widetilde{Q}\right)  \rightarrow
H_{k-1}\left(  \widetilde{Q},\partial\widetilde{U}_{i}\right)  \rightarrow
\cdots
\]
Since $H_{k}\left(  \widetilde{U}_{i},\partial\widetilde{U}_{i}\right)  $ and
$H_{k-1}\left(  \widetilde{Q},\partial\widetilde{U}_{i}\right)  $ are trivial
for all $k\leq n-3$, $H_{k}\left(  \widetilde{U}_{i},\widetilde{Q}\right)  $
also vanishes for $k\leq n-3$. If $k=n-2$, the surjectivity of $H_{n-2}\left(
\widetilde{Q},\partial\widetilde{U}_{i}\right)  \rightarrow H_{n-2}\left(
\widetilde{U}_{i},\partial\widetilde{U}_{i}\right)  $ together with the
triviality of $H_{n-3}\left(  \widetilde{Q},\partial\widetilde{U}_{i}\right)
$ implies the triviality of $H_{n-2}\left(  \widetilde{U}_{i},\widetilde
{Q}\right)  $. But the above $\pi_{1}$--isomorphisms imply that $\left(
\widetilde{V}_{i},\partial\widetilde{V}_{i}\right)  $ is the preimage of
$\left(  V_{i},\partial V_{i}\right)  $ under the covering projection
$p\co \left(  \widetilde{U}_{i},\partial\widetilde{U}_{i}\right)  \rightarrow
\left(  U_{i},\partial U_{i}\right)  $. Thus we may excise the interior of
$\widetilde{Q}$ from $\widetilde{U}_{i}$ to show that $H_{k}\left(
\widetilde{V}_{i},\partial\widetilde{V}_{i}\right)  $ vanishes for all $k\leq
n-2$.

\begin{remark}
With a few minor refinements, the above argument can be carried\break out when $n=5$
provided the Whitney Lemma is valid in the $4$--manifold\break $\partial_{+}\left(
S\cup h_{1}^{n-3}\cup\cdots\cup h_{t}^{n-3}\right)  $. This explains \fullref{r2}.
\end{remark}

\bibliographystyle{gtart}
\bibliography{link}

\end{document}